\documentclass[10pt]{article}
\usepackage{amsmath}
\usepackage{amssymb}
\usepackage{euscript}
\usepackage[latin1]{inputenc}

\normalbaselines
\newcommand{\qed}{\sqcap\kern-6.8pt\sqcup}

\newenvironment{proof}
      {\par\noindent{\it Proof\/: }\nopagebreak\normalsize}%
                                  {\linebreak[2]\hspace*{\fill}$\qed$\ifdim\lastskip<10pt
       \removelastskip \penalty-200  \vskip10pt  \fi}

\font\frten=eufm10 at 10pt
\font\freight=eufm10
\font\frsix=eufm8
\newfam\frfam\textfont\frfam=\frten
\scriptfont\frfam=\freight
\scriptscriptfont\frfam=\frsix

\newcommand{\CC}{\mathbb{C}}
\newcommand{\PP}{\mathbb{P}}
\newcommand{\RR}{\mathbb{R}}
\newcommand{\NN}{\mathbb{N}}
\newcommand{\ZZ}{\mathbb{Z}}

\newtheorem{thm}{Theorem}[section]
\newtheorem{prop}[thm]{Proposition}
\newtheorem{cor}[thm]{Corollary}
\newtheorem{lem}[thm]{Lemma}
\newtheorem{defn}[thm]{Definition}

\newtheorem{rema}[thm]{Remark}

\def \Clif{{\rm Clif}}
\def \dim{{\rm dim}}

\def \L{{\cal O}}

\def \Div{{\rm Div}}
\def \Gal{{\rm Gal}}

\def \deg{{\rm deg}}

\begin{document}

\title{Very special divisors on real algebraic curves}

\author{Jean-Philippe Monnier\\
       {\small D\'epartement de Math\'ematiques, Universit\'e d'Angers,}\\
{\small 2, Bd. Lavoisier, 49045 Angers cedex 01, France}\\
{\small e-mail: monnier@tonton.univ-angers.fr}}
\date{}
\maketitle

\begin{abstract}
{We study special linear systems called ``very special''
 whose dimension does not satisfy a
 Clifford type inequality given by Huisman. We classify all these very special
 linear systems when they are compounded of an
 involution. Examples of very special linear systems that are simple are also
 given.}
\end{abstract}

\section{Introduction and preliminaries}

\subsection{Introduction}

In this note, a real algebraic curve $X$ is a smooth proper geometrically
integral scheme over $\RR$ of dimension $1$.
A closed point $P$ of $X$ 
will be called a real point if the residue field at $P$ is $\RR$, and 
a non-real point if the residue field at $P$ is $\CC$.
The set of real points $X(\RR)$ of $X$ decomposes into 
finitely many connected components, whose number will be denoted by $s$.
By Harnack's Theorem (\cite[Th. 11.6.2 p. 245]{BCR}) we know that $s\leq g+1$,
where $g$ is the genus of $X$. 
A curve with $g+1-k$ real connected components
is called an $(M-k)$-curve. Another topological invariant associated to
$X$ is $a(X)$, the number of connected components of $X(\CC)\setminus X(\RR)$
counted modulo $2$. The pair $(s,a(X))$ is referred 
to as the topological type of $X$. A theorem of Klein asserts that there 
exists real curves of genus $g$ with topological type $(s,a)$ if and
only if the integers $g,s$ and $a$ obey the following restrictions:

\begin{prop}\cite{K}\label{klein}
\begin{description}
\item[1)] If $a(X)=0$, then $1\leq s\leq g+1$ and $s= g+1\mod 2$.
\item[2)] If $a(X)=1$, then $0\leq s\leq g$.
\end{description}
\end{prop}

The group $\Div(X)$ of divisors on $X$ 
is the free abelian group generated by the 
closed points of $X$. If $D$ is a divisor on $X$, we will denote by
$\L (D)$
its associated invertible sheaf. The dimension of the space 
of global sections of this sheaf will be denoted
by $\ell (D)$. 
Since a principal divisor has an even degree on each connected 
component of $X(\RR)$ (e.g. \cite{G-H} Lem. 4.1),
the number $\delta(D)$ (resp. $\beta(D)$) of connected components $C$
of $X(\RR )$
such that the degree of the restriction of $D$ to $C$ is odd (resp even) 
is an invariant of the linear system $|D|$
associated to $D$. If $\ell (D)>0$, the dimension of the linear system
$|D|$ is $\dim\, |D|=\ell (D)-1$.
Let $K$ be the canonical divisor. If 
$\ell (K-D)=\dim\, H^1 (X,\L (D) )>0$, $D$ is said to be special. If not,
$D$ is said to be non-special. By Riemann-Roch, if $\deg(D)>2g-2$ then $D$
is non-special. Assume $D$ is effective of degree $d$. If
$D$ is non-special then the dimension of the linear system $|D|$ is given by 
Riemann-Roch. If $D$ is special, then the dimension of the linear system
$|D|$ satisfies 
$$\dim\, |D|\leq \frac{1}{2} d.$$
This is the well known Clifford inequality for complex curves that
works also for real curves. 

Huisman (\cite[Th. 3.2]{Hu})
has shown that: 
\begin{thm} \label{Huisman}
Assume $X$ is an $M$-curve or an $(M-1)$-curve.
Let $D\in \Div(X)$ be an effective and special divisor of degree $d$.
Then $$\dim\, |D|\leq \frac{1}{2} (d-\delta (D)).$$
\end{thm}

Huisman inequality
is not valid for all real curves and the author has obtained the
following theorem.
\begin{thm} \cite[Th. A]{Mo}
\label{cliffreelA}
Let $D$ be an effective and special divisor of degree $d$. Then either
$$\dim\, |D|\leq \frac{1}{2} (d-\delta (D))\eqno{\rm (Clif 1)}$$
or 
$$\dim\, |D|\leq \frac{1}{2} (d-\beta (D))\eqno{\rm (Clif 2)}$$
Moreover, $D$ satisfies the inequality $(\Clif \,1)$
if either $s\leq 1$ or $s\geq g$.
\end{thm}

In this note we are interested in special divisors that do not satisfy
the inequality ${\rm (Clif 1)}$ given by Huisman.
\begin{defn} Let $D$ be an effective and special divisor of degree
  $d$. We say that $D$ is a very special divisor (or $|D|$ is a very
  special linear system)
  if $D$ does not
  satisfy the inequality (Clif 1) i.e. $\dim\, |D|>\frac{1}{2}
  (d-\delta (D))$.
\end{defn}

In the previous cited paper, the author has obtained a result in this direction.
\begin{thm} \cite[Th. 2.18]{Mo}
\label{siegalite}
Let $D$ be a very special and effective divisor of degree $d$ 
on a real curve $X$ such that (Clif 2) is an equality i.e.
$$r=\dim\, |D|= \frac{1}{2} (d-\beta (D))>\frac{1}{2} (d-\delta (D))$$
then $X$ is an hyperelliptic curve with 
$\delta( g_2^1 )=2$ and $|D|=rg_2^1 $ with $r$ odd.
\end{thm}

The author wishes to express his thanks to D. Naie and M. Coppens for
several helpful comments concerning the paper.

\subsection{ Preliminaries}

We recall here some classical concepts and notation we will
be using throughout the paper.

Let $X$ be a real curve. We will denote by $X_{\CC}$ the base extension of $X$
to $\CC$. 
The group $\Div(X_{\CC})$ of divisors on $X_{\CC}$
is the free abelian group on the 
closed points of $X_{\CC}$. The Galois group $\Gal(\CC /\RR)$
acts on the complex variety $X_{\CC}$ and also on $\Div(X_{\CC})$. 
We will always indicate this action by a bar. Identifying $\Div(X)$
and $\Div(X_{\CC})^{\Gal(\CC /\RR)}$,
if $P$ is a non-real point of $X$ then 
$P=Q+\bar{Q}$ 
with $Q$ a closed point of $X_{\CC}$. 

Let $D\in \Div(X)$ be a divisor with the property that $\L (D)$ has at
least one nonzero global section. The linear system $|D|$ is called
base point free if $\ell(D-P)\not=\ell(D)$ for all closed points
$P$ of $X$. If not, we may write $|D|=E+|D'|$ with $E$ a non zero
effective divisor called the base divisor of $|D|$, and with $|D'|$
base point free. A closed point
$P$ of $X$ is called a base point of $|D|$ if $P$ belongs to the support of
the base divisor of $|D|$. We note that $$\dim\, |D|=\, \dim\, |D' |.$$

As usual, a $g_d^r$ is an $r$-dimensional complete linear system
of degree $d$ on $X$. Let $|D|$ be a base point free $g_d^r$ on $X$.
The linear system $|D|$ defines a morphism $\varphi :\, X\rightarrow
\PP_{\RR}^r$ onto a non-degenerate (but maybe singular) curve in
$\PP_{\RR}^r$. If $\varphi$ is birational (resp. an isomorphism) onto
$\varphi (X)$, the $g_d^r$ (or $D$) is called simple (resp. very
ample). Let $X'$ be the normalization of $\varphi(X)$, and assume $D$ is not
simple i.e. $|D-P|$ has a base point for any closed point $P$ of $X$. Thus,
the induced morphism 
$\varphi :X\rightarrow X'$ is a non-trivial covering
map of degree $k\geq 2$. In particular, there exists $D'\in \Div(X')$ such that
$|D'|$ is a $g_{\frac{d}{k}}^{r}$ and such that $D=\varphi^* (D')$, i.e.
$|D|$ is induced by $X'$. If $g'$ denote the genus of $X'$,
$|D|$ is classically called compounded of an involution of order $k$
and genus $g'$. In the case $g'>0$, we speak of an irrational involution 
on $X$.

The reader is referred to \cite{ACGH} and
\cite{Ha} for more details on special divisors. Concerning real curves,
the reader may consult \cite{G-H}.
For $a\in\RR$ we denote by $[a]$ the integral part of $a$, i.e. 
the biggest integer $\leq a$.

\section{Non-simple very special divisors}

We first characterize the very special pencils.
\begin{prop}
\label{vspecpencil}
Let $D$ be a very special divisor of degree $d>0$ such that 
$\dim\, |D| =1$. Then $D=P_1 +\ldots + P_{s}$ with $P_1 ,\ldots ,P_{s} $
some real points of $X$ 
such that no two of them belong to the same connected component
of $X(\RR )$ i.e. $d=\delta (D)=s$. Moreover $D$ is base point
free.
\end{prop}

\begin{proof}
Since $D$ is special, we may assume that $D$ is
effective. Consequently, $d\geq \delta (D)$ and since 
$\dim\, |D| =1>\frac{1}{2} (d-\delta (D))$ we have $d=\delta (D)$
and $\dim\, |D| =1=\frac{1}{2} (d-\delta (D))+1$.

Since
$d=\delta (D)$, $D=P_1 +\ldots + P_{d} $ with $P_1 ,\ldots ,P_{d} $
some real points of $X$ 
such that no two of them belong to the same connected component
of $X(\RR )$.

Assume $d <s$. Choose a real point $P$ in one
of the $s-d $ real connected components that do not contain 
any of the points $P_1 ,\ldots ,P_{d} $. Since $\ell(D)=2$ then $\L
(D-P)$ has a nonzero global
section and $D-P$ should be linearly equivalent to an effective divisor
$D'$ of degree $d -1$ satisfying $ \delta (D') =d +1$.
This is impossible.

So $d=s$ and suppose $|D|$ is not base point free. If
$|D|$ has a real base point $P$, then $\dim\, |D-P|=1$ and 
$\deg (D-P)=\delta (D-P)=s-1$, contradicting the case $d<s$. 
If $|D|$ has a non-real
base point $Q$, then $\ell (D-Q)>0$ and $D-Q$ 
is linearly equivalent to an effective divisor
$D'$ of degree $s -2$ satisfying $ \delta (D') =s$,
which is again impossible.
\end{proof}

From the previous proposition, we get the the following corollary:
\begin{cor}
\label{1parcomp}
Let $D$ be a divisor of degree $d>0$ such that $D=P_1 +\ldots + P_{d}
$ with $P_1 ,\ldots ,P_{d} $
real points of $X$ 
such that no two of them belong to the same connected component
of $X(\RR )$.
Then $\dim\, |D| =0$ if $d <s$ and $\dim\, |D| \leq 1$ if $d =s$.
\end{cor}

The following lemma will allow us to restrict the study to base point free
linear systems.
\begin{lem} 
\label{reductbasepointfree}
Let $D\in \Div(X)$ be an effective divisor of degree $d$. Let $E$ be
the base divisor of $|D|$.
Let $|D'|=|D-E|$ be the degree $d'$ base point free part of $|D|$. 
If
$\dim\, |D' |\leq \frac{1}{2} (d'-\delta (D'))+k$
for a non-negative integer $k$, then
$\dim\, |D |\leq \frac{1}{2} (d-\delta(D))+k .$
\end{lem}

\begin{proof}
Write $D=D'+E$ where $E$ is the base divisor of $|D|$.
Assume $D'\in \Div(X)$ is an effective divisor of degree $d'$
satisfying $$\dim\, |D' |\leq \frac{1}{2} (d'-\delta (D'))+k$$
for a non-negative integer $k$.
Since  $\dim\, |D |=\dim\, |D'|$ and $E$ is effective, we have
$\delta
(D'+E)\leq\delta (D')+\deg(E) $.
Then $\dim\, |D |=\dim\, |D' +E | \leq \frac{1}{2} (d'-\delta (D'))+k
\leq \frac{1}{2} (\deg(D')+\deg(E) -\delta (D')-\deg(E))+k 
\leq \frac{1}{2} (\deg(D'+E) -\delta (D'+E))+k $ proving the lemma.
\end{proof}

Let $D$ be a special divisor.
Recall that $\delta (D)=\delta (K-D)$ and that $\beta (D)=\beta (K-D).$
The next lemma will allow us to study very special divisors 
of degree $\leq g-1$. 

\begin{lem}
\label{residuel}
Let $D\in \Div(X)$ be an effective and special divisor of degree $d$.
If
$\dim\, |D |= \frac{1}{2} (d-\delta (D))+k$
for a positive integer $k$, then
$\dim\, |K-D |= \frac{1}{2} (\deg(K-D)-\delta(K-D))+k .$
\end{lem}

\begin{proof} It is a straightforward calculation using Riemann-Roch.
\end{proof}

We can establish one of the main result of the paper.
\begin{thm}
\label{theorem1}
Let $D$ be a non-simple very special divisor of degree $d$.
Then $$\delta(D)=s$$ and $$\dim\, |D|= \frac{1}{2} (d-\delta
(D))+1.$$ Moreover $D$ is base point free.
\end{thm}

\begin{proof}
We prove, by induction on $\dim\, |D|$, the theorem for a base point
free non-simple very special
divisor. 

Let $D$ be a base point free non-simple divisors of degree $d$ such
that $\dim\, |D|=r>\frac{1}{2} (d-\delta(D))$. Since $D$ is special,
we may assume $D$ effective.

If $r=1$, Proposition \ref{vspecpencil} gives the result.

Assume $r>1$. Consider the map 
$\varphi:X\rightarrow \PP_{\RR}^r$
associated to $|D|$. Let $X'$ be the normalization of $\varphi (X)$. Then 
the induced morphism $\varphi :X\rightarrow X'$ is a non-trivial covering
map of degree $t\geq 2$ and there is $D'\in \Div(X')$ such that
$|D'|$ is a $g_{\frac{d}{t}}^{r}$ and such that $D=\varphi^* (D')$.

Assume $\delta(D)<s$. Let $P$ be a point of a connected component of
$X(\RR)$ where the degree of the restriction of $D$ is even. Since
$r\geq 1$, we may assume $D-P$ effective. Since $P$ is real,
$P'=\varphi (P)$ is real. Let $D_1=D-\varphi^*(P')$ and denote by $d_1=d-t$
its degree. Then $D_1$ is non-simple and effective since
$D_1=\varphi^*(D'-P')$. Moreover $\dim\, |D_1|=\dim\, |D-P|=r-1$.
Since $\dim\, |D-P|=r-1=\frac{1}{2} (d-\delta(D))+1-1=
\frac{1}{2} (d-1-(\delta(D)+1))+1=\frac{1}{2}
(\deg(D-P)-\delta(D-P))+1$, we see that $D-P$ is a very special
divisor. By Lemma \ref{reductbasepointfree}, $D_1$ is also very
special. Since $D$ is base point free, then $D'$ is base point free.
Choosing $P$ such that $D'-P'$ is base point free then $D_1$ is base
point free since $\dim\, |D_1=\varphi^*(D'-P')|=\dim\, |D'-P'|$.
By induction, $\delta(D_1)=s$ and 
$\dim\, |D_1|=r-1=\frac{1}{2} (d_1-s)+1$ i.e.
\begin{equation}
\label{equ1}
r=\frac{1}{2} (d_1-s)+2.
\end{equation}
Remark that $d\geq d_1+2$ ($\varphi$ is non-trivial) and
$\delta(D)\leq s-1$. 
If 
$\delta(D)= s-1$ then $d\geq d_1+3$, since $\delta(D_1)=s $ and
$\deg(\varphi^* (P'))\geq 2$. Hence we get
$r>\frac{1}{2} (d-\delta(D))\geq 
\frac{1}{2} (d_1+3-\delta(D))= \frac{1}{2} (d_1+3-s+1)=\frac{1}{2}
(d_1-s)+2$ contradicting (\ref{equ1}).
If $\delta(D)< s-1$,
we get $r>\frac{1}{2} (d-\delta(D))\geq 
\frac{1}{2} (d_1+2-\delta(D))\geq \frac{1}{2} (d_1+2-s+2)=\frac{1}{2}
(d_1-s)+2$ contradicting (\ref{equ1}).

We have just proved that $\delta(D)=s$. Now assume $r\geq
\frac{1}{2} (d-s)+2$. Let $P$ be a real point. Let $D_1$ be the
divisor of degree $d_1$ constructed as in the above proof
that $\delta(D)=s$.
Since $\dim\, |D-P|=r-1\geq\frac{1}{2} (d-s)+2-1=
\frac{1}{2} (d-1-(\delta(D)-1))+1=\frac{1}{2}
(\deg(D-P)-\delta(D-P))+1$, we see that $D-P$ is a very special
divisor. By Lemma \ref{reductbasepointfree}, $D_1$ is also very
special. For a general choice of $P$, $D_1$ is also base point free.
By induction, $\delta(D_1)=s$ and 
$\dim\, |D_1|=r-1=\frac{1}{2} (d_1-s)+1$ since $D_1$ is non-simple and
base point free. Since $d\geq d_1+2$ then $r=\frac{1}{2} (d_1-s)+2=
\frac{1}{2} (d_1+2-s)+1\leq \frac{1}{2} (d-\delta(D))+1$, impossible.

We have proved the theorem in the case of base point free divisors.
Let $D$ be a non-simple very special divisor. If $D$ has base point,
with the previous notation, it means that $D'$ has base point.
Write $D=D_2+E$ where $E$ is the base divisor of $|D|$.
Since $D=\varphi^* (D')$ and $\dim\, |D|=\dim\, |D'|$, we have 
$E=\varphi^*(E')$ where $E'$ is the base divisor of $|D'|$, it means
that $D_2$ is also non-simple.
By Lemma \ref{reductbasepointfree} and the proof for base point free divisors,
$D_2\in \Div(X)$ is an effective divisor of degree $d_2$
satisfying $\dim\, |D_2 |=r= \frac{1}{2} (d_2-\delta (D_2))+1$
and $\delta(D_2)=s$. Let $e$ denote the degree of $E$. Since
$\dim\, |D_2 |= \frac{1}{2} (d_2-\delta (D_2))+1$ and since $D$ is
very special, we have $r=  \frac{1}{2} (d-\delta(D))+1$
by Lemma 
\ref{reductbasepointfree}. But $r= \frac{1}{2} (d_2-s)+1$, hence
$d-d_2=e=\delta(D)-s\leq 0$. Consequently, $e=0$ and $\delta(D)=s$,
i.e. $D$ is base point free.
\end{proof}

\section{Very special nets}

A net is a linear system of dimension $2$ i.e. a
$g_d^2$.

We recall some classical definitions concerning real curves in
projective spaces.
Let $X\subseteq \PP_{\RR}^{r}$, $r\geq 2$, 
be a real curve. $X$ is
non-degenerate if $X$ is not contained in any hyperplane of
$\PP_{\RR}^r$. In what follows, $X$ is supposed to be non-degenerate.
Let $C$ be a connected component of $X(\RR )$.
The component $C$ is called a pseudo-line if the canonical class of $C$
is non-trivial in $H_1 (\PP_{\RR}^r (\RR ),\ZZ /2)$.
Equivalently, $C$ is
a pseudo-line
if and only if for each real hyperplane $H$,
$H(\RR )$ intersects $C$ in an odd number of points, 
when counted with multiplicities (see \cite{Hu}).

Before looking at very special nets, we need some lemmas concerning
morphisms between real curves and very special divisors.
\begin{lem} 
\label{basic1}
Let $\varphi :X\rightarrow X'$ be a covering map of degree $t$
between two real curves $X$ and $X'$.
\begin{description}
\item[({\cal{i}})] If $P'$ is a real point of $X'$ then 
  $\varphi^{-1} (P')$ can contain real and non-real points.
\item[({\cal{ii}})] If $Q'$ is a non-real point of $X'$ then 
  $\varphi^{-1} (Q')$ is totally non-real.
\item[({\cal{iii}})] The image by $\varphi$ of a connected component $C$ of
  $X(\RR)$ is either a connected component of $X' (\RR)$, or 
  a compact connected semi-algebraic subset of a connected
  component $C'$ of $X'(\RR )$ corresponding topologically to a closed
  interval of $C'$.
\end{description}
\end{lem}

The proof of the lemma is trivial.

\begin{lem}
\label{surjreel}
Let $\varphi :X\rightarrow X'$ be a covering map of degree $t$
between two real curves $X$ and $X'$. Let $P$ be a real point of
$X(\RR )$ contained in a connected component $C$ of $X(\RR
)$. Let $C'$ denote the connected component of $X'(\RR)$ containing
$P'=\varphi(P)$. If $\deg(\varphi^{*} (P')\cap C)$ is odd then $\varphi(C)=C'$.
\end{lem}

\begin{proof}
Suppose $\varphi (C)$ is not a
connected component of $X'(\RR)$. Then 
$\varphi(C)$, corresponds topologically to a closed
interval of a connected component of $X'(\RR)$. Let $P_1'$ be
one of the two end-points of this interval. Let $P_1\in \varphi^{-1} (P_1')\cap
C$. Then the ramification index $e_{P_1}$ is even since $C$ is clearly
on one side of the fiber $\varphi^{-1} (P_1')$. Hence  
$\deg(\varphi^{*} (P_1')\cap C)$ is even. It is impossible since 
$\deg(\varphi^{*} (P_1')\cap C)= \deg(\varphi^{*} (P')\cap C)\mod 2$.
\end{proof}

The following lemma is a generalization of a lemma due to Huisman \cite{Hu}.
\begin{lem}
\label{huisman2}
Let $D\in  \Div(X)$ be a divisor of degree $d$ such that $\ell(D)>0$. 
Assume that $d+\delta(D)< 2s +2k$ with
$k\in \NN$. Then $$\dim\, |D|\leq \frac{1}{2} (d-\delta (D))+k.$$
\end{lem}

\begin{proof}
We proceed by induction on $k$.
The case $k=0$ is given by Lemma
\cite{Hu}.
So, assume that $k>0$ and that
$d+\delta(D)<2s+2k$. Since $\ell(D)>0$, we may assume that $D$ is effective.
If $d+\delta(D)<2s+2k-2$, the proof is done by the induction hypothesis.
Since $d=\delta (D)\,\mod\, 2$, we assume that $d+\delta(D)=2s+2k-2$.
Let $Q$ be a non-real point. Since 
$\deg(D-Q)+\delta(D-Q)<2s+2k-2$ and $\delta(D-Q)=\delta (D)$, if $\ell
(D-Q)>0$ then, using the induction hypothesis,
$\dim\, |D-Q|\leq \frac{1}{2} (\deg(D-Q)-\delta (D-Q))+k-1=\frac{1}{2}
(d-2-\delta (D))+k-1$. 
Hence $\dim\, |D|\leq \dim\, |D-Q|+2\leq \frac{1}{2} (d-\delta (D))+k$.
If $\ell (D-Q)=0$ then $\dim\, |D|\leq 1 \leq \frac{1}{2} (d-\delta
(D))+k$, 
since we have $k>0$ and $d\geq\delta (D)$ ($D$ is effective).
\end{proof}

The following proposition shows that the excess
can be bounded in terms of $r$ for linear systems of 
dimension $r$ which do not satisfy (Clif 1).
\begin{prop}
\label{dimension}
Let $D$ be an effective and special divisor of degree $d$ on a
real curve $X$. Assume that $r=\dim\, |D|= \frac{1}{2} (d-\delta (D))+k+1.$
Then $k\leq [\frac{r-1}{2}]$.
\end{prop}

\begin{proof}
Assume $r=2n-\varepsilon$ with $\varepsilon\in\{ 0,1\}$.
We proceed by induction on $n$.

If $r=1$, we get $2=d-\delta (D)+2k+2$. Since $d\geq \delta (D)$, 
we have $k=0$ and $d=\delta (D)$.

If $r=2$, we get $4=d-\delta (D)+2k+2$ i.e. $2=d-\delta (D)+2k$.
If $k\geq 1$, we must have $d=\delta (D)$, hence $d+\delta (D)\leq 2s$.
But Lemma \ref{huisman2} says that $k=0$, which is impossible.
So $k=0$ and $d=\delta (D)+2$.

Assume $n>1$ and $k>0$. Choose two real points $P_1 ,P_2$ in the same real 
connected component such that $\dim\, |D-P_1 -P_2 |=r-2$.
Then $r-2=2(n-1)-\varepsilon =\frac{1}{2} ((d-2)-\delta (D-P_1 -P_2 ))+(k-1)+1.$
By the induction hypothesis, $(k-1)\leq \bigg[\frac{(r-2)-1}{2}\bigg]$ i.e.
$$k\leq \bigg[\frac{r-1}{2}\bigg].$$
\end{proof}

The following theorem is the main result of this section.

\begin{thm}
\label{thm2} 
If $D$ is a very special divisor then $\dim\, |D |\not=2$.
\end{thm}

\begin{proof}
Let $D$ be a very special divisor such that $\dim\, |D |=2$. By
Lemma \ref{reductbasepointfree}, we can assume $D$ is base point free.
By Proposition \ref{dimension} we get
$2= \frac{1}{2} (d-\delta(D))+1$, i.e. $d=\delta(D)+2$.
By Lemma \ref{huisman2}, we have $d+\delta(D)\geq 2s$. Hence
$\delta(D)\geq s-1$ and we have two possibilities: either $d=s+2$ and
$\delta(D)=s$, or $d=s+1$ and
$\delta(D)=s-1$.\\

First, assume $D$ is simple.
In this case, $X$ is mapped birationally by $\varphi$---associated to
$|D|$---onto 
a curve of degree $d$ in $\PP_{\RR}^2$. By the genus formula,
$$g=\frac{1}{2}(d-1)(d-2)-\mu,$$ with $\mu$ the multiplicity of the
singular locus of $\varphi (X)$. 
If $d=s+2$, we know that $\varphi (X)$ has
exactly $s$ pseudo-lines. Since any two distinct pseudo-lines of
$\varphi (X)$ intersect each other, we have $\mu\geq
\frac{1}{2}(s-1)(s)$. By the genus formula, $g\leq
\frac{1}{2}(s+1)(s)-\frac{1}{2}(s-1)(s)=s$. Hence $X$ is an $M$-curve
or an $(M-1)$-curve and Theorem \ref{Huisman} gives a contradiction.
If $d=s+1$, we similarly find $g\leq s-1$. Hence $X$ is an $M$-curve and Theorem
\ref{Huisman} leads again to a contradiction.\\

Second, assume $D$ is not simple.
By
Theorem \ref{theorem1}, $\delta(D)=s$ and
$2=\frac{1}{2} (d-s)+1$ i.e. $d=s+2$.
Consider the map 
$f:X\rightarrow \PP_{\RR}^r$
associated to $|D|$. Let $X'$ be the normalization of $f (X)$. Then 
the induced morphism $\varphi :X\rightarrow X'$ is a non-trivial covering
map of degree $t\geq 2$ and there is $D'\in \Div(X')$ such that
$|D'|$ is a $g_{\frac{d}{t}}^{r}$ and such that $D=\varphi^* (D')$.
Let $C_1 ,\ldots,C_s$ denote the connected components of $X(\RR
)$. Since $\dim\, |D |=2$, we may assume $D=P_1+\ldots+P_s+R_1+T_1$
with $P_i,R_i,T_i\in C_i$. Since the support of $D$ is totally real,
Lemma \ref{basic1} implies that the support of $D'$ is totally real.
By Lemma \ref{surjreel}, $\varphi (C_i)$ is a connected component of
$X'(\RR)$ for $i\geq 2$.

Suppose $t=2$. The points $P_1$, $R_1$, $T_1$ verify 
that no two of them belong to the same fiber over the points of $D'$
because if it is not the case the degree $\mod 2$ of the restriction
of the fiber to $C_1$ is not constant. Hence we may assume that we
have $3$ points 
$P_1'$, $R_1'$, $T_1'$ of the support of $D'$ contained in the same
connected component $C'$ of $X'(\RR)$ such that
$\varphi^{*}(P_1')=P_1+P_2$, $\varphi^{*}(R_1')=R_1+P_3$, 
$\varphi^{*}(T_1')=T_1+P_4$. Then by Lemma
\ref{surjreel},
$\varphi(C_1)=\varphi(C_2)=\varphi(C_3)=\varphi(C_4)=C'$ which is clearly impossible since it would imply that 
$P_2,P_3,P_4$ belong to the same connected component.

Suppose $t\geq 3$. Arguing similarly to the case $t=2$, we conclude
that $\varphi(P_1)=\varphi(R_1)=\varphi(T_1)$. Hence the degree
of the restriction to
a connected component of $X(\RR)$ of a fiber over a point from the
support of $D'$ is either empty or of odd degree. By
Lemma \ref{surjreel}, $\varphi (X(\RR))$ is a union of connected
components of $X'(\RR)$. Moreover, since
$\varphi(P_1)=\varphi(R_1)=\varphi(T_1)$ and since 
no two points among $P_1 ,\ldots ,P_{s} $
belong to the same connected component
of $X(\RR )$, we get that $\deg (D')=\delta (D')$. The support of
$D'$ consists on a single point in each connected component of 
$\varphi (X(\RR))$. Corollary \ref{1parcomp} says that $\dim\, |D'
|=\dim\, |D
|\leq 1$, which is again impossible.
\end{proof}

\section{Construction of the non-simple very special divisors}

We recall the definition of a $g'$-hyperelliptic curve 
(see \cite[p. 249]{FK}).
A curve $X$ is called $g'$-hyperelliptic if there exists $\varphi
:X\rightarrow X'$ a non-trivial covering map of degree $2$ such that the
genus of $X'$ is $g'$.
Classical hyperelliptic curves correspond to $0$-hyperelliptic curves.\\

In this section, we prove that the non-simple very special linear
systems of dimension $r>1$ are lying on some ``very special''
$g'$-hyperelliptic curves and its converse.
\begin{thm}
\label{construction}
Let $D$ be a non-simple very special divisor of degree $d$ such that
$\dim\, |D|=r>1$. Let $\varphi :X\rightarrow X'$ be the non-trivial covering
map of degree $t\geq 2$
induced by $|D|$ such that there is $D'\in \Div(X')$ of degree $d'$ such that
$|D'|$ is a $g_{\frac{d}{t}}^{r}$ and such that $D=\varphi^*
(D')$. Let $g'$ denote the genus of $X'$ and let $s'$ be the number of
connected components of $X'(\RR)$.
Then
\begin{description}
\item[({\cal{i}})] $D$ is base point free, $r=\frac{1}{2}(d-s)+1$ and
  $\delta(D)=s$;
\item[({\cal{ii}})] $t=2$ i.e. $X$ is a $g'$-hyperelliptic curve;
\item[({\cal{iii}})] $s$ is even, $s'=\frac{s}{2}$, $\varphi
  (X(\RR))=X'(\RR)$, the inverse image by
  $\varphi$ of
  each connected component of $X'(\RR)$ is a disjoint union of $2$
  connected components of $X (\RR)$;
\item[({\cal{iv}})] $r$ is odd and $\delta(D')=s'$;
\item[({\cal{v}})] $D'$ is a base point free non-special divisor and
  $X'$ is an $M$-curve;
\item[({\cal{vi}})] $D'$ is linearly
  equivalent to an effective divisor
  $P_1'+\ldots+P_{s'}'+R_{1,2}'+\dots+R_{1,r}'$ 
  with $P_i',R_{i,j}'\in C_i'$ such that $\dim\,|P_1'+\ldots+
  P_{s'}'|=1$ and such that $R_{1,2}',\ldots,R_{1,r}'$ are general in
  $C_1'$, 
  where 
  $C_1',\ldots, C_{s'}'=C_{\frac{s}{2}}'$ denote
  the connected components of $X'(\RR)$;
\item[({\cal{vii}})] $a(X)=0$ and $g$ is odd;
\item[({\cal{viii}})] there is a very special pencil on $X$: 
$|\varphi^*(P_1'+\ldots+P_{s'}')|$.
\end{description}
\end{thm}

\begin{proof}
We keep the notation and the hypotheses of the theorem.
Theorem \ref{theorem1} gives statement ({\cal{i}}).

Assume $t\geq 3$. 
Let $Q'$ be a non-real point of $X'$. Since $r>2$ (by Theorem
\ref{thm2}), we may assume $Q'$ 
in the support of $D'$. By Lemma \ref{basic1},
  $\delta(\varphi^*(Q'))=0$. Let $D_1=D-\varphi^* (Q')$ and let $d_1$
  denote its degree. We have $\dim\,|D_1|=\dim\,|D'-Q'|=r-2=
\frac{1}{2}(d-s)+1-2=\frac{1}{2}(d-4-s)+1>\frac{1}{2}(d-2t-s)+1=
\frac{1}{2}(\deg(D_1)-\delta(D_1))+1$. Theorem \ref{theorem1} says
that this is not possible and then statement ({\cal{ii}}) is proved.

Since $t=2$ then $d=2d'$. Since $r=\frac{1}{2}(d-s)+1$ we see that $s$
is even.
Suppose $\varphi (X(\RR))\not= X'(\RR)$. Using Lemma
\ref{surjreel}, there exists a real point $P'$ of $X'(\RR)$ such that
$\varphi^{-1}(P')$ is non-real. Let $D_1=D-\varphi^* ((r-2)P')$ and let $d_1$
  denote its degree. Choosing $P'$ sufficiently general, $\dim\,|D_1|=
\dim\,|D'-(r-2)P'|=2=\frac{1}{2}(d-s)+1-(r-2)=
\frac{1}{2}(d-2r+4-s)+1=
\frac{1}{2}(\deg(D_1)-\delta(D_1))+1$. Theorem \ref{thm2} says
that this is not possible, hence  $\varphi (X(\RR))= X'(\RR)$.

If $C'$ is a connected component of $X'(\RR)$, we have two
possibilities for the inverse image by
  $\varphi$ of $C'$: either $\varphi^{-1}(C')$ is a disjoint union of $2$
  connected components of $X (\RR)$, or $\varphi^{-1}(C')$ is a
  connected component of $X(\RR)$. In the second case, choosing 
  a real point $P'$ of $C'$, we can prove as above that 
  $D_1=D-\varphi^* ((r-2)P')$ is a very special divisor such that
  $\dim\,|D_1|=2$, contradicting Theorem \ref{thm2}.
Hence the inverse image by
  $\varphi$ of
  each connected component of $X'(\RR)$ is a disjoint union of $2$
  connected components of $X (\RR)$. Consequently $s'=\frac{s}{2}$ and
  $\delta(D)=2\delta(D')$. Since $\delta(D)=s$, we have $\delta(D')=s'$.

We claim that $r$ is odd. Indeed, if $r$ is even, choosing 
  a non-real point $Q'$ of $X$ then we can prove as above that 
  $D_1=D-\varphi^* (\frac{r-2}{2}Q')$ is a very special divisor such that
  $\dim\,|D_1|=2$, again impossible by Theorem \ref{thm2} establishing
  the claim.

So $r$ is odd. We can choose general points $R_{1,2}',\dots,R_{1,r}'$
in $C_1'$ such
that $\ell (D'-(R_{1,2}'+\dots+R_{1,r}'))=2$. 
Since $\deg (D'-(R_{1,2}'+\dots+R_{1,r}'))=r-1+\frac{s}{2}-(r-1)=s'$,
we may assume that there are real
points $P_1',\ldots, P_{s'}'$ such that $P_i'\in C_i'$,
$D'=P_1'+\ldots+P_{s'}'+R_{1,2}'+\dots+R_{1,r}'$ and $\dim\,|P_1'+\ldots+
P_{s'}'|=1$. Moreover, $|\varphi^* (P_1'+\ldots+P_{s'}')|$ is a very
special pencil. We have proved the statements ({\cal{vi}}) and ({\cal{viii}}).

The divisor $D'$ is non-special. If $D'$ were special, 
$D_1'=P_1'+\ldots+P_{s'}'+R_{1,2}'$
would be special as a subdivisor of $D'$. Moreover
$\dim\,|D_1'|=2$ by the above construction. But
$\dim\,|D_1'|=\frac{1}{2}(\deg(D_1')-\delta(D_1'))+1$, hence $D_1'$
would be
a very special divisor such that $\dim\,|D_1'|=2$, a contradiction.

Since $D'$ is non-special,
$\dim\,|D'|=r=\frac{1}{2}d-\frac{1}{2}s+1=d'-s'+1=d'-g'$ by Riemann-Roch. Hence
$s'=g'+1$ and $X'$ is an $M$-curve.

Since $X'$ is an $M$-curve, $a(X')=0$ (see Proposition
\ref{klein}). Since $\varphi^{-1} (X'(\RR))=X (\RR)$ then $\varphi
(X(\CC)\setminus X(\RR))=X'(\CC)\setminus X'(\RR)$. If
$X(\CC)\setminus X(\RR)$ is connected then also $X'(\CC)\setminus X'(\RR)$,
impossible. Hence $a(X)=0$. Since $a(X)=0$ and $s$ is even then $g$ is
odd by Proposition \ref{klein}.
\end{proof}

\begin{cor}
\label{separante}
If $X$ has a non-simple very special divisor then $a(X)=0$.
\end{cor}

\begin{proof}
Let $D$ be a non-simple very special divisor of degree $d$ such that
$\dim\, |D|=r$.
If $r>1$, Theorem \ref{construction} gives the result.
If $r=1$, let $\varphi :X\rightarrow\PP_{\RR}^1$ be the morphism
induced by $|D|$. By Lemma \ref{vspecpencil}, we have $\varphi
^{-1}(\PP_{\RR}^1(\RR))=X(\RR)$. Since $a(\PP_{\RR}^1)=0$, we easily get
$a(X)=0$.
\end{proof}

\begin{rema}
\label{exitrigo}
{\rm By Theorem \ref{construction}, if there is a non-simple very
  special divisor $D$ on $X$ such that $\dim\,|D|>1$ then there is a
  very special pencil on $X$. The converse is not true. For example,
  let $X$ be real trigonal curve such that $\delta (g_3^1)=3$. By
  \cite[p. 179]{G-H}, such a trigonal curve exists. The $g_3^1$ is
  very special and we get $s=3$ by Proposition
  \ref{vspecpencil}. Since $s$ is odd, Theorem \ref{construction}
  says that there is not a non-simple very
  special divisor $D$ on $X$ such that $\dim\,|D|>1$.}
\end{rema}

\begin{rema}
{\rm In the situation of Theorem \ref{construction}, if in addition
  the genus of $X'$ is $0$, then $X$ is an hyperelliptic curve with
  $\delta(g_2^1)=2$ and $|D|=rg_2^1$ with $r$ odd. Such 
  hyperelliptic curves exist (see \cite[Rem. 2.11]{Mo}) for any odd
  genus $g\geq 3$
  and such very
  special divisors have already been studied in this particular case,
  they correspond to the extremal cases in \cite[Th. 2.18]{Mo}.}
\end{rema}

In the rest of the section, we prove the converse of 
Theorem \ref{construction}.
We first put a name on the curves appearing in Theorem
\ref{construction}.
\begin{defn}
\label{defveryspec}
A curve $X$ of genus $g$ is a very special
$g'$-hyperelliptic curve if $g$ is odd and if  
there exists $\varphi :X\rightarrow X'$ a non-trivial covering
map of degree $2$ such that $X'$ is an $M$-curve of genus $g'$ and 
the inverse image by $\varphi$ of any connected component of $X'(\RR)$
is the union of two connected components of $X(\RR)$ i.e. $X$ satisfies all the
topological properties of Theorem \ref{construction}.
\end{defn}

The existence of very special $g'$-hyperelliptic curves of odd genus
$g$ is proved in \cite{BEG}.
\begin{prop} 
\label{exiveryspec}
\cite{BEG}\\
Let $g'$, $g$ be natural numbers, $g\geq 2$. For each odd $g$, and
each $g'$ verifying $2g'+2\leq g+1$, there exists a very special
$g'$-hyperelliptic curve of genus $g$.
\end{prop}

\begin{proof}
The existence of very special $g'$-hyperelliptic curves of odd genus
$g$ is not explicitely written in \cite{BEG}. We explain here how to
deduce from \cite[Th. 2]{BEG} and from the proof of
\cite[Thm. 2]{BEG} the existence of these curves.

Let $g'$, $g$ be natural numbers, $g\geq 2$. By \cite[Th. 2]{BEG},
if $2g'+2\leq g+1$, there exists a $g'$-hyperelliptic curve $X$ of
genus $g$ whose real part has $s=2g'+2$ connected components. Furthemore
$a(X)=0$.

Now assume $g$ odd. We have to look at the proof of
\cite[Th. 2]{BEG} in order to see if wether or not the curve $X$ built in this
proof is very special. The curve $X$ is very special if
the number of connected components of $X'(\RR)$ is $g'+1$ and if 
$\varphi :X\rightarrow X'$ has no real branch point. These two
conditions are satisfied (see the case (b1) p. 280 of the proof
of \cite[Th. 2]{BEG} and see \cite{BCG} for an explanation of the
notations).
\end{proof}

Since we have the existence of very special $g'$-hyperelliptic curves
of odd genus $g$ verifying $2g'+2\leq g+1$, we prove now the converse
of Theorem \ref{construction} on these very special curves.

\begin{prop}
\label{invconstruction}
Let $X$ be a very special $g'$-hyperelliptic curves
of odd genus $g$ such that $2g'+2\leq g$.
Let $\varphi :X\rightarrow X'$ denote the
corresponding non-trivial covering
map of degree $2$.
Let $C_1',\ldots, C_{g'+1}'$ denote
the connected components of $X'(\RR)$. 
Let $P_1',\ldots,P_{g'+1}$ be some real points of $X'(\RR)$ such that
$P_i'\in C_i'$.
If $r\geq 2$, let $R_{1,2}',\dots,R_{1,r}'$ be general points in $C_1'$.
We set $D'_1=P_1'+\ldots+P_{g'+1}'$,
$D'_r=P_1'+\ldots+P_{g'+1}'+R_{1,2}'+\dots+R_{1,r}'$ for $r\geq 2$,
and $D_r=\varphi^*
(D_r')$.\\
For any odd $r$
such that $$1\leq r\leq g-2g'-1$$
choosing the general points $R_{1,2}',\dots,R_{1,r}'$
such that the two real points of $\varphi^{-1}(R_{1,j}')$ are not base
points of $K-D_{j-1}$, $2\leq j\leq r$,
then $D_r$ is a non-simple very
special 
divisor
such that $\dim |D_r|=r$.
\end{prop}

\begin{proof}
If $g'=0$ i.e. $X$ is hyperelliptic and $\varphi$ is the hyperelliptic
map, then \cite[Prop. 2.10]{Mo} gives the result for any choice 
of $R_{1,2}',\dots,R_{1,r}'$.\\

For the rest of the proof, we assume $g'>0$.
\\
{\it Claim 1:} $|D_1|$ is a very special pencil 
since $2g'+2\leq g$.\\
We consider the linear system $|D_1|$.
By Riemann-Roch, we have $\dim |D_1'|\geq 1$ hence $\dim |D_1|\geq 1$.
In fact $\dim |D_1|= 1$ by Corollary \ref{1parcomp} and since the
support of $D_1$ consists of exactly one point in each connected
component of $X(\RR)$. If $D_1$ is non-special, then 
$$\dim |D_1|=1=s-g=2g'+2-g$$ by Riemann-Roch. We get a
contradiction with the hypothesis $2g'+2\leq g$.
\\
{\it Claim 2:} If $2\leq r\leq g-2g'-1$ then $D_r$ is special
for any choice of $R_{1,2}',\dots,R_{1,r}'$.\\
By Riemann-Roch $\dim |D_r|\geq \dim |D_r'|\geq r$.
Hence $\dim |D_r|=r+l$ with an integer $l\geq 0$.
Assume $D_r$ non-special. We get $r+l=2g'+2r-g$ by
Riemann-Roch,
i.e. $r=g-2g'+l$, a contradiction.
\\
Now we prove by induction on $r\geq 1$ that $\dim |D_r|=r$.\\
If $r=1$, Claim 1 gives the result.
Now suppose $2\leq r+1 \leq g-2g'-1$.
By Claim 2, $D_{r+1}$ is
special. Notice that $D_{r}$ is also
special since $D_{r}$ is an effective subdivisor of $D_{r+1}$.
By the induction hypothesis, we get $\dim |D_r|=r$.
Let $\varphi^* (R_{1,r+1}')=R_1+R_2$ then
$D_{r+1}=D_{r}+R_1+R_2$. 
Assume $\dim |D_{r+1}|>r+1$ then clearly
$\dim |D_{r+1}|=r+2$. Moreover 
$\dim |D_r+R_1|=r+1$ and $D_r+R_1$ is special since it is an effective 
subdivisor of $D_{r+1}$. Let $K$ denote the canonical divisor of $X$.
Since $\dim |D_r+R_1|=\dim |D_r|+1$, then $R_1$ is a base point of
$|K-D_r|$ contradicting the general choice of $R_{1,r+1}'$ and
finishing
the induction.\\
\\
If $$1\leq r\leq g-2g'-1,$$ then, by  Claim 1, Claim 2 and the
induction argument with the
points $R_{1,2}',\dots,R_{1,r}'$ general, $D_r$ is
special and non-simple such that $\dim |D_r|=r$. If $r$ is
odd, it is easy to see that $\dim |D_r|=\frac{1}{2}(\deg
(D_r)-\delta(D_r))+1$ and the proof is done.
\end{proof}

Coppens has made the following remark concerning a consequence of
Proposition \ref{invconstruction}.
\begin{cor}
\label{remacoppens}
For any odd $g\geq 2$ and any even $0<d<g$, there exists a real
curve $X$ of genus $g$ with a very special pencil of degree $d$.
\end{cor}

One interesting question of Coppens is wether the result of Corollary
\ref{remacoppens} is also valid for even genus.
Using a result in \cite{BEGM}, we give a partial answer to that
question.

\begin{prop}
\label{partialcoppens}
Let $g\geq 2$. For any $0<d<g$ such that there exists an integer
$k\geq 2$ with $g=(k-1)(d-1)$, there exists a real
curve $X$ of genus $g$ with a very special pencil of degree $d$.
\end{prop}

\begin{proof}
If $$g-1-(k-1)d+k=0$$
i.e. if $g=(k-1)(d-1)$,
by \cite{BEGM}[Prop. 2], there exists a real
curve $X$ of genus $g$ with $d$ real connected components with a cyclic
morphism $f:X\rightarrow \PP_{\RR}^1$ having only
ramification points of index $d$ over $k$ non-real
points of $\PP_{\RR}^1$.
Clearly, $f$ corresponds to a very special pencil of degree $d$.
\end{proof}

A consequence of the previous Proposition is a different proof 
of a Gross and Harris \cite[p. 179]{G-H} result
mentioned in Remark \ref{exitrigo} of the previous
section, concerning the existence of trigonal curves with a very
special pencil.
\begin{cor}
\label{exitrigocyclic}
For every even $g\geq 4$, there exists a real trigonal curves of genus
$g$ with a $g_3^1$ very special.
\end{cor}

\section{Simple very special divisors and very special curves in some
  projective spaces}

In this section, we prove the existence of simple very special
divisors.\\

\begin{prop}
\label{simpletrigo} Let $X$ be a real trigonal curve and let $D$ be a
divisor on $X$ such that:
\begin{description}
\item[({\cal{i}})]
 $|D|=g_3^1$,
\item[({\cal{ii}})]
$\delta (g_3^1)=3$,
\item[({\cal{iii}})]
$g\geq 5$.
\end{description}
Then, the base point free part $D'$ of $K-D$ is a simple very special
divisor with $\dim\,|D'|=g-3$.
Moreover $g$ is even and $s=3$.
\end{prop}

\begin{proof}
As we have already noticed in Corollary \ref{exitrigocyclic} 
such a trigonal curve with $\delta (g_3^1)=3$ exists.
Since $g\geq 5$, the $g_3^1$ is unique.
By Proposition \ref{vspecpencil}, we get $s=3$.
Since $a(X)=0$, we see that $g$ is even.
If $|D|=g_3^1$ then $D$ is a non-simple very special divisor.
Let $D'$ be the base point free part of $K-D$.
Since $\ell(K-D')\geq \ell(D)=2$, $D'$ is a special divisor.
By Lemmas \ref{residuel} and \ref{reductbasepointfree}, $D'$ is also
very special. By Riemann-Roch
$$\dim\,|D'|=\dim\,|K-D|=(2g-5)-g+2=g-3\geq 2.$$
Since $s$ is odd, Theorem \ref{construction} forces $D'$ to be simple.
\end{proof}

Let $X\subseteq \PP_{\RR}^{r}$, $r\geq 2$, 
be a smooth real curve. We assume, in what follows, that $X$ is non-degenerate.
We say that $X$ is special (resp. very special) if the divisor
associated to the sheaf of hyperplane sections $\L_{X} (1)$
is special (resp. very special).

\begin{cor} For every odd $r\geq 3$, there exists a very special curve in 
$ \PP_{\RR}^{r}$.
\end{cor}

\begin{proof} Let $r\geq 3$ be an odd integer.
Let $X$ be a real trigonal curve and let $D$ be a
divisor on $X$ such that:
\begin{description}
\item[({\cal{i}})]
 $|D|=g_3^1$,
\item[({\cal{ii}})]
$\delta (g_3^1)=3$,
\item[({\cal{iii}})]
$g=r+3$.
\end{description}
By Proposition \ref{simpletrigo} the base point free part $D'$ of 
$K-D$ is a simple very special
divisor with $\dim\,|D'|=r$. Hence $\varphi_{|D'|}(X)$ is a 
very special curve in 
$ \PP_{\RR}^{r}$.
\end{proof}

\end{document}